\RequirePackage{snapshot}
\documentclass[10pt]{amsart}

\usepackage{amsmath,amssymb,amsfonts,epsfig,color,euscript}

\newcommand{\ig}{\includegraphics}
\newcommand{\rb}{\raisebox}
\newcommand\risS[6]{\rb{#1pt}[#5pt][#6pt]{\begin{picture}(#4,15)(0,0)
  \put(0,0){\ig[width=#4pt]{#2.eps}} #3
     \end{picture}}}

\newtheorem{thm}{Theorem}[section]
\newtheorem{defn}[thm]{Definition}
\newtheorem{exa}[thm]{Example}

\newtheorem{lemma}[thm]{Lemma}

\newtheorem{cor}[thm]{Corollary}

\theoremstyle{remark}
\newtheorem*{rem}{Remark}

\begin{document}

\title[On a conjecture of Gross, Mansour and Tucker]
  {On a conjecture of Gross, Mansour and Tucker}
\author[Sergei Chmutov, Fabien Vignes-Tourneret]{
       Sergei Chmutov, Fabien Vignes-Tourneret}

\date{\today}
\address{(S.Ch.) The Ohio State University at Mansfield,
1760 University Drive Mansfield, OH 44906, USA.
{\tt chmutov.1@osu.edu}}
\address{(F.V.-T.) Institut Camille Jordan, Universit\'e de Lyon, CNRS UMR 5208 ;
Bât. Braconnier, 43 bd du 11 novembre 1918, F-69622 Villeurbanne Cedex, France.
{\tt vignes@math.univ-lyon1.fr}}


\begin{abstract}
Partial duality is a duality of ribbon graphs relative to a subset of their edges generalizing the classical Euler–Poincar\'e duality. This operation often changes the genus.
Recently J.~L.~Gross, T.~Mansour, and T.~W.~Tucker formulated a conjecture that for any ribbon graph different from plane trees and their partial duals, there is a subset of edges partial duality relative to which does change the genus. A family of counterexamples was found by Qi Yan and Xian'an Jin. In this note we prove that essentially these are the only counterexamples.
\end{abstract}

\maketitle

\section{Ribbon graphs and partial duality} \label{s:rg-pd}
A {\it\bfseries ribbon graph} topologically can be regarded as a regular neighborhood of the graph cellularly embedded into a closed surface. Then the neighborhoods of vertices and edges represent {\it\bfseries vertex-discs} and 
{\it\bfseries edge-ribbons}. A component of the complement of the ribbon graph is called a {\it\bfseries face-disc}. The formal definition is the following:
\begin{defn}\label{def:rb}{\rm
A {\it\bfseries ribbon graph} $G$ is a surface with boundary represented as the union of two sets of closed topological discs called
{\it\bfseries vertex-discs} $V(G)$ and {\it\bfseries edge-ribbons} $E(G)$, satisfying the following
conditions:
\begin{itemize}
\item[$\bullet$] the vertex-discs and edge-ribbons intersect by disjoint line segments;
\item[$\bullet$] each such line segment lies on the boundary of precisely one vertex-disc and precisely one edge-ribbon;
\item[$\bullet$] every edge-ribbon contains exactly two such line segments.
\end{itemize} 
}\end{defn}

\noindent
The last condition allows us to view each edge-ribbon as a topological rectangle 
(ribbon) attached to vertex-discs by two opposite sides. The face-discs are assumed to be discs which are glued to each of the boundary component of the ribbon graph to get a closed surface $\Sigma_G$.

In this paper we are dealing with orientable ribbon graphs only. Also this paper deals a lot with one-vertex ribbon graphs which we will encode by 
{\it\bfseries chord diagrams} following the common practice in knot theory
\cite[Sec.4.8]{CDM}. Namely we attach a vertex-disc to the
circle of a chord diagram and thicken the chords to narrow edge-ribbons.
Here is an example
$$\risS{-9}{cd34}{}{24}{0}{0}\qquad \risS{-1}{totorl}{}{20}{0}{0}\qquad
\risS{-25}{basket}{}{50}{35}{25}\ =\ 
\risS{-20}{gr34}{}{70}{0}{0}\ =\ G,\label{basket}
\hspace{1cm}
\Sigma_G\ =\ \risS{-25}{rg-torus}{}{100}{0}{0}
$$

\noindent
This ribbon graph has two boundary components. The corresponding closed surface 
$\Sigma_G$ has thus two face-discs and its Euler characteristic is equal to
$\chi(\Sigma_G)=1-3+2=0$. This is a torus and the genus of $G$ is equal to 1.

{\it\bfseries Partial duality} of ribbon graphs was introduced in \cite{Ch} under the name of generalized duality. Then it was further studied and developed in papers \cite{VT,Mo1,Mo2,EMM12,BBC12,Mo3,HM13,CVT, EZ17, GMT,YaJi20}.
We refer to \cite{EMM} for an excellent exposition.

\begin{defn}\label{def:pd}{\rm For a ribbon graph $G$ and a subset of its edge-ribbons $A\subseteq E(G)$, the {\it\bfseries partial dual} $G^A$ of $G$ relative to $A$ is a ribbon graph constructed in the following way. The vertex-discs of $G^A$ are bounded by connected components of the boundary of the spanning subgraph of $G$ containing all the vertices of $G$ and only the edges from $A$. The edge-ribbons of
$E(G)\setminus A$ are attached to these new vertices exactly at the same places as in $G$. The edge-ribbons from $A$ become parts of the new vertex-discs now. 
For an $e\in A$ we take a copy of $e$, $e'$, and attach it to the new vertex-discs in the following way. The rectangle representing $e$ intersects with vertex-discs of $G$ by a pair of opposite sides. But it intersects the boundary of the spanning subgraph, that is the new vertex-discs, along the arcs of the other pair of its opposite sides. We attach $e'$ to these arcs by
this second pair. The copies of the first pair of sides in $e'$ becomes the arcs of the boundary of $G^A$.
}\end{defn}
The following figure illustrates this definition for the case of an edge $e$ connecting different vertices (the boxes with dashed arcs mean that there might be other edges attached to these vertices).
$$G =\ 
  \risS{-8}{dus1}{\put(34,-2){\mbox{$e$}}}{70}{0}{0}\quad 
   \risS{-4}{totor}{}{30}{0}{0}\quad 
  \risS{-35}{dus3}{\put(-4,13){\mbox{$e'$}}}{110}{45}{40}\ =\ 
  \risS{-35}{dus4}{\put(47,65){\mbox{$e'$}}}{70}{10}{0}\ =\ G^{\{e\}}
$$
Partial duality relative to a loop acts on this figure backwards, from right to left.
In the example from previous page the partial duality relative to
middle loop gives:
$$G\ = \  \risS{-20}{gr34}{\put(50,49){\mbox{$e$}}}{70}{40}{25}\sim
           \risS{-9}{cd34}{}{24}{0}{0}\qquad 
\risS{-1}{totor}{}{20}{0}{0}\qquad \risS{-9}{cd34-pd}{}{24}{0}{0}\sim
\risS{-18}{gr34-pd}{\put(40,17){\mbox{$e'$}}}{70}{0}{0}\ =\ G^{\{e\}}\ ,
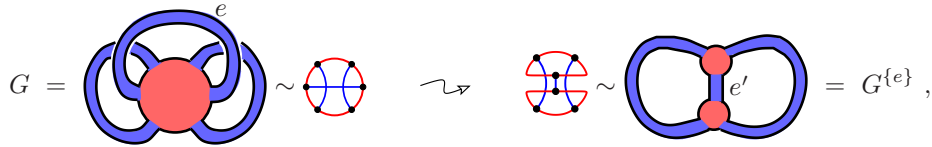
$$
a graph cellularly embedded into a sphere.

Partial duality relative to a set of edges can be done step by step one edge at a time. Partial duality relative to the whole set of edges $E(G)$ is the  classical Euler–Poincar\'e duality. Partial duality relative to a spanning tree produces a one-vertex ribbon graph which can be encoded by a chord diagram.

We see that partial duality relative to single edge $e\in E(G)$ changes the structure of vertices attached to the edge-ribbon $e$. Namely, if there were two separate vertex-discs attached to $e$ in $G$, then they merge together to single vertex-disc attached to the corresponding edge-ribbon $e'$ in $G^{\{e\}}$, and 
vice versa. A careful analysis of the orientation of the boundary of the ribbon graphs following the arrows on the above picture shows that a similar rearrangement happened with the face-discs. We formulate it as a lemma.

\begin{lemma}
If the edge-ribbon  $e\in E(G)$ was attached to two distinct vertex- (resp. face-)  discs of $G$, then in the ribbon graph  $G^{\{e\}}$ these discs will merge together into a single vertex- (resp. face-) disc of $G^{\{e\}}$ attached to the corresponding edge-ribbon $e'\in E(G^{\{e\}})$. Backwards, if $e\in E(G)$ was attached to a single vertex- (resp. face-) disc of $G$, then in the ribbon graph  
$G^{\{e\}}$ it will be split into two vertex- (resp. face-) discs of $G^{\{e\}}$ attached to the corresponding edge-ribbon $e'\in E(G^{\{e\}})$. 
\end{lemma}

According to the lemma there are four types of edges depending on the
attachment to the vertex- and face-discs. The impact of partial
duality relative to a single edge $e$ on the genus $g(G)$ depends on the type. These four types correspond to the upper four rows (oriented case) of Table 1.1 of \cite{GMT}. We will follow the notation of \cite{GMT} denoting each type with two letters $p$ or $u$, where $p$ stands a {\it proper} edge connecting distinct vertices and $u$ stands for an 
{\it untwisted loop}. The second letter indicates the similar properties of the corresponding edge of the dual graph.

{\bf Type $pp$.} {\it\bfseries $e$ is attached to different vertices $v_1\neq v_2$ and to two different faces $f_1\neq f_2$.} In this case the partial duality relative to $e$ merges the vertices and the faces and creates a loop $e'$ attached to the resulting vertex and resulting face. So the genus $g(G)$ will be increased by 1, and the corresponding edge-ribbon $e'$ in $G^{\{e\}}$ falls into Type $uu$.

{\bf Type $uu$.} {\it\bfseries $e$ is a loop attached to a single vertex $v_1=v_2$ and to a single face $f_1=f_2$.} The partial duality splits the vertex and the face into two vertices $v_1$, $v_2$ and two faces $f_1$, $f_2$. The edge $e$ returns to Type $pp$, decreasing the genus $g(G)$ by 1.

{\bf Type $pu$.} {\it\bfseries $e$ is attached to two different
  vertices $v_1\neq v_2$ but to a single face $f_1=f_2$.} The partial duality relative to $e$ merges the vertices but split the face into two faces $f_1$, $f_2$. The genus $g(G)$ stays the same, and the  corresponding edge-ribbon $e'$ in $G^{\{e\}}$ falls into Type $up$.

{\bf Type $up$.} {\it\bfseries $e$ is a loop attached to a single vertex $v_1=v_2$ and to two different faces $f_1\neq f_2$.} The partial duality splits the vertex and merge the faces. The edge $e$ returns to Type $pu$ and the genus $g(G)$ stays the same.

The next figures illustrate partial duality for these four types.
\begin{figure}[!htp]
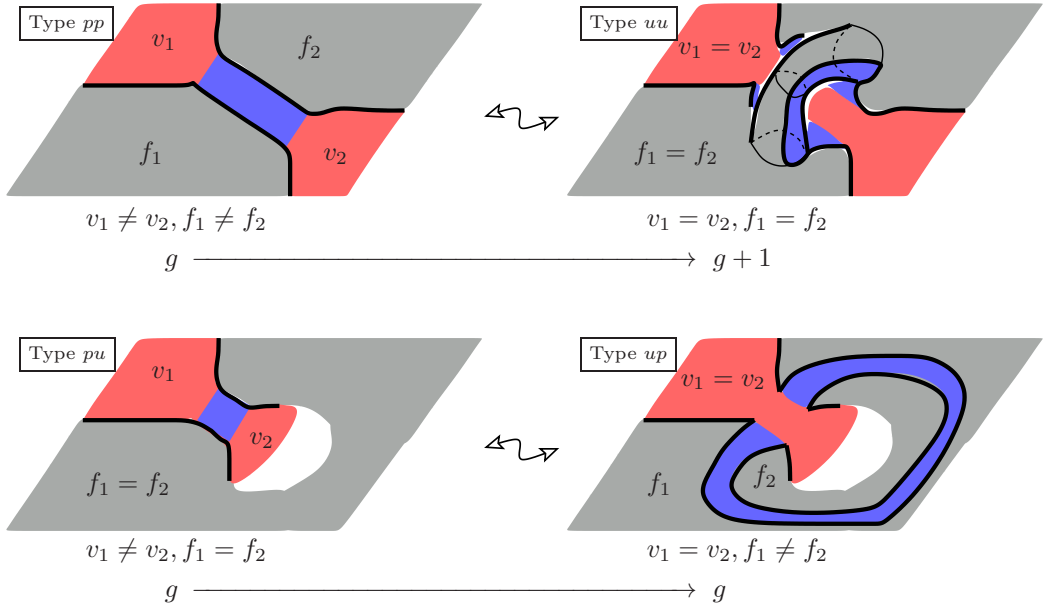

$$\qquad\risS{-30}{pd-1a}{\put(5,65){\fbox{\scriptsize Type $pp$}}
  \put(55,58){\mbox{$v_1$}}\put(120,15){\mbox{$v_2$}}
  \put(50,15){\mbox{$f_1$}}\put(110,55){\mbox{$f_2$}}
  \put(30,-10){\mbox{$v_1\neq v_2, f_1\neq f_2$}}
  \put(60,-25){\mbox{$g$}}
  \put(70,-25){\mbox{$\xrightarrow{\hspace*{6.5cm}}$}}}{180}{50}{0}
\risS{-4}{totorl}{}{30}{0}{0}\
\risS{-30}{pd-2a}{\put(5,65){\fbox{\scriptsize Type $uu$}}
  \put(42,55){\mbox{$v_1=v_2$}}\put(25,15){\mbox{$f_1=f_2$}}
  \put(30,-10){\mbox{$v_1=v_2, f_1=f_2$}}
  \put(55,-25){\mbox{$g+1$}}}{180}{35}{80}\qquad
$$

$$\qquad\risS{-30}{pd-3a}{\put(5,65){\fbox{\scriptsize Type $pu$}}
  \put(55,58){\mbox{$v_1$}}\put(92,33){\mbox{$v_2$}}
  \put(30,15){\mbox{$f_1=f_2$}}
  \put(30,-10){\mbox{$v_1\neq v_2, f_1= f_2$}}
  \put(60,-25){\mbox{$g$}}
  \put(70,-25){\mbox{$\xrightarrow{\hspace*{6.5cm}}$}}}{180}{35}{0}
\risS{-4}{totorl}{}{30}{0}{0}\
\risS{-30}{pd-4a}{\put(5,65){\fbox{\scriptsize Type $up$}}
  \put(43,55){\mbox{$v_1=v_2$}}\put(70,18){\mbox{$f_2$}}
  \put(30,15){\mbox{$f_1$}}
  \put(30,-10){\mbox{$v_1=v_2, f_1\neq f_2$}}
  \put(55,-25){\mbox{$g$}}}{180}{35}{70}\qquad
$$
\caption{Types of edges in partial duality}
\label{fig-edgetypes}
\end{figure}

\section{The Gross--Mansour--Tucker conjecture} \label{s:GMT-conj}

In \cite{GMT} J.~L.~Gross, T.~Mansour, and T.~W.~Tucker introduced the 
{\it\bfseries partial-dual orientable genus polynomial} ${\ }^\partial\Gamma_G(z)$ 
as a generating function of the numbers of partial duals of $G$ of the given genus:
$${\ }^\partial\Gamma_G(z):=\sum_{A\subseteq E(G)} z^{g(G^A)}$$
where $g(G^{A})$ stands for the orientable genus of $G^{A}$, namely
half its Euler genus.

The {\it\bfseries ribbon-join operation} $G_1\vee G_2$ on ribbon graphs $G_1$ and $G_2$ was introduced by I.~Moffatt \cite{Mo3}. To obtain $G_1\vee G_2$ glue together
a vertex-disc of $G_1$ and a vertex-disc of $G_2$ along some arcs on
their boundaries between consecutive ends of the end of edge ribbons
to get a new vertex disc of $G_1\vee G_2$. For one vertex ribbon
graphs this operation is the standard multiplication of corresponding
chord diagrams, see \cite[Sec.\@ 4.4.2]{CDM}. Obviously the result is not unique, it depends on choice of vertices and arcs. But,
since the partial dualities of factors $G_1$ and $G_2$ are
independent, namely
\begin{equation*}
  \forall E_{1}\subset E(G_{1}), E_{2}\subset E(G_{2}), \ (G_{1}\vee
  G_{2})^{E_{1}\cup E_{2}}=G_{1}^{E_{1}}\vee G_{2}^{E_{2}},
\end{equation*}
the polynomial ${\ }^\partial\Gamma_G(z)$ is multiplicative with respect to this operation, ${\ }^\partial\Gamma_{G_1\vee G_2}(z)= {\ }^\partial\Gamma_{G_1}(z) \cdot 
{\ }^\partial\Gamma_{G_2}(z)$.\\

J.~L.~Gross, T.~Mansour, and T.~W.~Tucker formulated the following conjecture.

{\bf GMT-conjecture.} \cite[Conjecture 3.1]{GMT} {\it  There is no orientable ribbon graph having a non-constant partial-dual genus polynomial
with only one non-zero coefficient.}

\medskip
They noted that if $G$ is a tree, which implies that its genus is 0, all partial duals of $G$ also will be of genus 0. That means for a tree with $n$ edges,
${\ }^\partial\Gamma_G(z)=2^n$. This explains the exclusion of constants from the conjecture.

Q.~Yan, X.~Jin \cite{YaJi20} found a family of counterexamples to the GMT-conjecture.
They used the one-vertex ribbon graphs $B_n$ with $n$ edge-ribbons every two of which are interlaced along the boundary circle of the vertex. These are the corresponding first chord diagrams.
$$B_1=\risS{-9}{cd-K1}{}{24}{15}{15},\qquad
B_2=\risS{-9}{cd-K2}{}{24}{0}{0},\qquad
B_3=\risS{-9}{cd-K3}{}{24}{0}{0},\qquad
B_4=\risS{-9}{cd-K4}{}{24}{0}{0},\qquad
B_5=\risS{-9}{cd-K5}{}{24}{0}{0},\qquad\dots
$$

An {\it\bfseries intersection graph} of a chord diagram
\cite[Sec.4.8]{CDM} is a graph whose vertices correspond to the chords
and two vertices are connected by an edge if the corresponding chords intersect. The ribbon graphs $B_n$ can be characterized as one-vertex ribbon graphs for which chord diagrams have the complete graph $K_n$ as the intersection graph.
 
It was proved in \cite[Theorem 22]{YaJi20} that for odd values of $n$, 
${\ }^\partial\Gamma_{B_n}(z)=2^nz^{\frac{n-1}2}$, thus providing the counterexample to the GMT-conjecture. 

\begin{exa}{\rm To demonstrate that the even values of $n$ do not
    work, we need to check the attachment of edge-ribbon to the face-discs. For that we double each chord modeling the boundary sides of the corresponding rectangles and look for the boundary components we obtain. For example,
$$B_2=\risS{-9}{cd-K2}{}{24}{15}{15}\qquad \risS{-1}{totor}{}{20}{0}{0}\qquad
\risS{-9}{cd-K2-dch}{}{24}{15}{15}\qquad \risS{-1}{totor}{}{20}{0}{0}\qquad
\risS{-9}{cd-K2-f}{}{24}{15}{15}
$$ 
We can see a single immersed curve representing the boundary of the ribbon graph
$B_2$. Each boundary component is a circle to which a face-disc is supposed to be glued in to form a closed surface. In this case we have only one face-disc. So each ribbon of $B_2$ is attached to the same face. Hence it is of Type $uu$, and the partial dual relative to each would change the genus. Consequently, the partial-dual genus polynomial would have more than one term. Meanwhile for odd values of $n$,
for example $n=3$, each ribbon would be attached to different boundary components:
$$B_3=\risS{-9}{cd-K3}{}{24}{15}{15}\qquad \risS{-1}{totor}{}{20}{0}{0}\qquad
\risS{-9}{cd-K3-dch}{}{24}{15}{15}\qquad \risS{-1}{totor}{}{20}{0}{0}\qquad
\risS{-9}{cd-K3-f}{}{24}{15}{15}
$$ 
We would have two immersed curves, so we would have corresponding two face-discs, and each edge ribbon would be attached to different face-discs. That means they are of Type $up$ and the partial duality relative to a single edge would not change the genus. To complete the calculation of \cite{YaJi20} we would need to check that
all edge-ribbons will stay withing Types $pu$ and $up$ after performing partial dualties.
}\end{exa}

By the multiplicativity property, the partial-dual genus polynomial of the 
ribbon-join of $B_n$'s with odd $n$'s also will be a monomial. In fact, even for trees, the total duality gives the one-vertex ribbon graph for which the chords of its chord diagram do not intersect each other. So it is a crossing-less matching of $2n$ points on a circle. As such it is a ribbon-join of $n$ copies of $B_1$.
Therefore, all the know examples of one-vertex ribbon graph with a monomial as the
partial-dual genus polynomial consist of ribbon-joins of various $B_n$'s with odd values of $n$. Our main theorem claims that these are the only counterexamples to the GMT-conjecture.

\section{Main result} \label{s:main-r}

We call a connected ribbon graph $G$ {\it\bfseries join-prime} if it
cannot be represented as the ribbon-join of two
graphs $G_{1},G_{2}$ with at least one 
edge-ribbon each: $G \neq G_1\vee G_2$.

\begin{thm}\label{thm-main} For any join-prime ribbon graph different from partial duals of 
$B_n$'s with odd values of $n$, there are partial duals of different genus.
\end{thm}

One can observe that there are exactly two graphs in the partial dual class of 
$B_n$. Namely $B_n$ itself and its partial dual relative to one edge-ribbon. Partial duality relative to two edges of $B_n$ will return to 
$B_n$ back. So, for any join-prime ribbon graph different from these two the
GMT-conjecture holds. Moreover, from the multiplicativity of ${\
}^\partial\Gamma$ with respect to the join operation, Theorem
\ref{thm-main} implies that the only orientable ribbon graphs having a
non-constant partial-dual genus polynomial with only one non-zero
coefficient are partial duals of joins of $B_{n}$'s with $n$ odd.

\begin{rem}
  The non-orientable counterpart of the Gross-Mansour-Tucker
  conjecture has been treated by Maya Thompson (Royal Holloway
  University of London). She proved that the only non-orientable
  ribbon graphs having a non-constant partial-dual genus polynomial with only one non-zero
  coefficient are partial duals of the one-vertex ribbon graph with
  one twisted edge. Moreover, a preprint by Q.\@ Yan and X.\@ Jin \cite{Yan2021as}
  posted recently contains in particular proofs of both the orientable and
  non-orientable cases.
\end{rem}

A one-vertex ribbon graph is join-prime if and only if the intersection graph of the corresponding chord diagram is connected.

Our proof of Theorem \ref{thm-main} is based on the following Lemma.
\begin{lemma} Let $G$ be a one-vertex join-prime ribbon graph and
  $e\in E(G)$. If for all subsets $A\subseteq E(G)$, $g(G)=g(G^A)$ then
\begin{enumerate}
\item $e$ is attached to different face-discs $f_1\neq f_2$. That is $e$ has to be of Type $up$.
\item Any edge-ribbon interlaced with $e$ is attached to the same face-discs 
$f_1$ and $f_2$. 
\item Any edge-ribbon not interlaced with $e$ is attached to a pair of face-discs different from $\{f_1,f_2\}$.
\end{enumerate}
\end{lemma}

\begin{proof}[\it\bfseries Proof of the Lemma.]
Since $g(G)=g(G^{\{e\}})$, the edge-ribbon $e$ has to be of Type
$up$. The same is true for any edge-ribbon of $G$. This proves
(1). Let $e_1$ be an edge-ribbon interlaced with $e$ and attached to
two different face-discs $f_3$ and $f_4$. Then on the partial dual
$G^{\{e\}}$ the corresponding edge $e_1$ will connect two different
vertices. If we want the genus to stay the same, $g(G)=g(G^{\{e,e_1\}})$,
the edge $e_1$ has to be of Type $pu$ in $G^{\{e\}}$. This means it has to be attached to the same face-disc of $G^{\{e\}}$. In other words the faces  $f_3$ and 
$f_4$ merge to each other in  $G^{\{e\}}$. But the only change in the structure of faces under the partial duality relative to $e$ is that the faces $f_1$ and $f_2$ merge to a single face. This implies $\{f_3,f_4\}=\{f_1,f_2\}$, which proves (2).
In a similar way, let $e_2$ be an edge-ribbon not interlaced with $e$
and attached to two different face-discs $f_5$ and $f_6$. Then, after
the partial duality relative to $e$, $e_{2}$ is still a loop attached
to one of the two vertices of $G^{\{e\}}$. So it has to be of Type $up$ in $G^{\{e\}}$ which means that
$f_5\neq f_6$ in $G^{\{e\}}$.  But this partial duality merges $f_1$ and $f_2$. So, 
$\{f_5,f_6\}\neq\{f_1,f_2\}$, which proves (3).
\end{proof}

\begin{proof}[\it\bfseries Proof of the Theorem.]
Since every connected ribbon graph has a one-vertex graph as partial dual, we assume
that $G$ is a one-vertex join-prime ribbon graph with $n$ edge-ribbons and such that 
$g(G)=g(G^A)$ for all subsets $A\subseteq E(G)$. Let $D_G$ be the
corresponding chord diagram. We also assume that the intersection
graph of $D_G$ is different from the complete graph $K_n$. Then there is a path of length 2 in the intersection graph. Let
$e_2,e_1,e_3$ be the vertices of the path. They correspond to chords of $D_G$ such that the chord $e_1$ intersects both $e_2$ and $e_3$, but $e_2$ and $e_3$ do not intersect each other. So, the chord diagram  $D_G$ looks like this
$$D_G=\risS{-9}{cd34-dot}{}{24}{15}{12},
$$ 
where the dashed circle means that there might be other chords with end
points on these arcs. In terms of the ribbon graph $G$, the
edge-ribbon $e_{1}$ is interlaced with both $e_2$ and $e_3$, but the edge-ribbons 
$e_2$ and $e_3$ do not interlace. Every 
edge-ribbon has to be of Type $up$, that is attached to two different face-discs.
Let $\{f_1,f_2\}$ (resp. $\{f_3,f_4\}$ and $\{f_5,f_6\}$) be the pair
of face-discs attached to $e_1$ (resp. $e_2$ and $e_3$). Then
according to case (2) of the Lemma, 
$\{f_3,f_4\}=\{f_1,f_2\}=\{f_5,f_6\}$. But this contradicts case (3) of the lemma according to which $\{f_3,f_4\}\neq\{f_5,f_6\}$. This means that any two vertices have to be actually adjacent. This implies that the intersection graph has to be the complete graph $K_n$ which proves the Theorem.
\end{proof}

T.~W.~Tucker indicated that the theorem can be reformulated in the following ways.
\begin{cor} The only join-prime ribbon graph $G$ where all edges have Type $up$ is a one-vertex ribbon graph with an odd number of edges where every loop interlaces with every other loop.
\end{cor}
Or in the dual form:
\begin{cor} The only join-prime ribbon graph where all edges have Type $pu$  is the dipole on an odd number of edges where the rotation at one vertex is the opposite of the rotation at the other vertex. (i.e. the Petrie dual of the standard dipole in the sphere).
\end{cor}

\section{Future research} \label{s:fut}
\begin{itemize}
\item The paper \cite{GMT} contains another very interesting Conjecture 5.1 that the sequence of coefficients of ${\ }^\partial\Gamma_G(z)$ is {\it log-concave}
as well as analogous conjectures for the partial-dual Euler-genus
polynomial in the non-orientable case.

\item The paper \cite{CVT} generalizes partial duality to hypermaps. It would be interesting to find (and prove) the analogous conjectures for hypermaps.

\item Ribbon graphs may be considered from the point of view of delta-matroids 
\cite{CMNR}. In this way the concepts of partial duality and genus can
be interpreted in terms of delta-matroids. It would be interesting to
know whether the partial-dual genus polynomial and the related conjectures would make sense for general delta-matroids.
\end{itemize}

\section*{Acknowledgment}
We are very grateful to J.~L.~Gross and T.~W.~Tucker for numerous
comments on the preliminary draft of this paper and to the anonymous
referees the suggestions of whom greatly improved the exposition.

\end{document}